\documentclass[CRMATH,Unicode,manuscript]{cedram}

\usepackage{epsfig}
\usepackage[hang]{caption2}
\usepackage[normalsize,it,hang]{subfigure}

\title{Flatness-based control revisited: The {\it HEOL} setting}  

\alttitle{Revoir la commande par platitude : la m\'{e}thode {\rm HEOL}}

\author{\firstname{C\'{e}dric} \lastname{Join}
\CDRorcid{0000-0001-6304-1262}}
\address{CRAN (CNRS, UMR 7039),
Universit\'{e} de Lorraine, Campus Aiguillettes, BP 70239, 54506 Vand{\oe}uvre-l\`{e}s-Nancy, France}
\email[C. Join]{cedric.join@univ-lorraine.fr, cedric.join@alien-sas.com}
 \addressSameAs{4}{ AL.I.E.N., 7 rue Maurice Barr\`{e}s, 54330 V\'{e}zelise, France 
 }
\author{\firstname{Emmanuel} \lastname{Delaleau}\CDRorcid{0000-0002-0257-568X}}
\address{ENI Brest, UMR CNRS 6027, IRDL, 29200 Brest, France}
\email[E. Delaleau]{delaleau@enib.fr}

\author{\firstname{Michel}
\lastname{Fliess}
\CDRorcid{0000-0003-3426-3637}}
\address{LIX (CNRS, UMR 7161),
\'{E}cole polytechnique, 91128 Palaiseau, France}
\email[M. Fliess]{michel.fliess@polytechnique.edu, michel.fliess@alien-sas.com, michel.fliess@swissknife.tech}
\address{AL.I.E.N., 7 rue Maurice Barr\`{e}s, 54330 V\'{e}zelise, France 
  }


\begin{abstract} 
We present the algebraic foundations of the HEOL setting, which combines flatness-based control and intelligent controllers, two advances in automatic control that have been proven in practice, including in industry. The result provides a solution to many pending questions on feedback loops concerning flatness-based control and model-free control (MFC). Elementary module theory, ordinary differential fields and the generalization of K\"{a}hler differentials to differential fields provide an intrinsic definition of the tangent linear system. The algebraic manipulations associated with the operational calculus lead to homeostat and intelligent controllers. They are illustrated via some computer simulations.
\end{abstract}

\begin{altabstract} 
On présente les fondations algébriques de la méthode HEOL qui combine commande par platitude et bouclage intelligent, c'est-à-dire deux avancées de l'automatique ayant fait leur preuve en pratique, y compris industrielle. On résoud ainsi plusieurs questions pendantes sur les bouclages à propos de la platitude et de la commande sans modèle. Théorie élémentaire des modules, corps différentiels ordinaires, et la généralisation à ces corps des différentielles de K\"{a}hler permettent une définition intrinsèque du système linéaire tangent. Les manipulations algébriques associées au calcul opérationnel conduisent à l'hom\'{e}ostat et aux correcteurs intelligents, illustrés par simulations numériques. \\\\
{\bf Mots-clefs.} Commande par platitude, commande sans modèle, contr\^{o}leurs intelligents, algèbre différentielle, théorie des modules.
\end{altabstract}

\begin{document}

\maketitle

\section{Introduction}
({\em Differentially}) {\em flat} systems \cite{flmr2,flmr4}, which were introduced more than thirty years ago \cite{flmr1}, have been undeniably influential, not only in control engineering (see, \textit{e.g.}, \cite{flmr2,flmr4}, the books \cite{levine,miunske,rigatos,rudolph,sira}, and references therein), but also in other fields such as pure physics (see, \textit{e.g.}, \cite{phys}). Their discovery is the result of years of research into feedback linearization of systems modeled by ordinary differential equations. Let us mention here a few important steps: static state-feedback linearization \cite{jakub}, \cite{hunt}, and dynamic feedback linearization \cite{charlet}.\footnote{See, \textit{e.g.}, \cite{clelland} for a survey on dynamic feedback linearization from the point of view of symmetry.} For flat systems, a special type of dynamic feedback, called \emph{endogenous}, is used. As well known (see, \textit{e.g.}, \cite{delaleau,hagen1,hagen2}) their implementation is difficult. Flat systems possess another characteristic that is as unexpected as it is essential: There exists a finite set $\{y_1, \dots, y_m\}$ of variables, called \emph{flat}, or \emph{linearizing}, \emph{output}, such that
\begin{itemize}
\item any system variable $z$ may be expressed as a \emph{differential function} of the component of the flat output and their derivatives up to some finite order, \textit{i.e.}, $$z = f(y_1, \dots, y_m, \dots, y_1^{(\nu_1)}, \dots, y_m^{(\nu_m)})$$
\item any component of the flat output may be expressed as a differential function of the system variables;
\item the components of the flat output are \emph{differentially independent}, \textit{i.e.}, they are not related by any differential relation.
\end{itemize}
Assigning time functions to $y_1, \dots, y_m$ yields time functions to any system variable  without any integration procedure:
\begin{itemize}
\item this is an \emph{open loop} or \emph{feedforward} control strategy; 
\item it provides a \emph{reference trajectory}.
\end{itemize}
This feature plays a prominent r\^{o}le in concrete applications of flatness-based control.
\begin{rema}
A control system is an \emph{underdetermined} system of ordinary differential equations, \textit{i.e.}, a system where the number of equations is less than the number of unknown variables. It is worthy of note to cite here a little-known paper \cite{hilbert} by Hilbert. He considers there a single differential equation with two unknown variables. In our terminology he asks when such a system is flat, \textit{i.e.}, when the unknowns may be obtained without any integration procedure. There is no hint of any linearization!  
\end{rema}

Any practitioner knows the difficulty if not the impossibility of writing down a ``good'' mathematical modeling in ``complex'' situations. It turns out however that ``(over)simplified'' modeling is quite often flat. They may be useful for deriving an open-loop reference trajectory. Several publications (see in chronological order \cite{villagra}, \cite{cancer}, \cite{covid}, \cite{mounier}, \cite{saar}, \cite{wang}) have successfully closed the loop via the \emph{intelligent} controllers associated to \emph{model-free control} (\emph{MFC}) \cite{mfc1,mfc2} to mitigate mismatches and disturbances. The {\em ultra-local model} \cite{mfc1,mfc2} reads in the case of the single control (resp. output) variable $u$ (resp. $y$)
\begin{equation}\label{ul}
    \frac{d^\nu}{dt^\nu} \Delta y = F + \alpha \Delta u
\end{equation}
where 
\begin{itemize}
    \item $\Delta y = y - y^\star$, $\Delta u = u - u^\star$ where $y^\star$ is a reference trajectory and $u^\star$ the corresponding reference control;
    \item the time-varying term $F$ outlines all the quantities which are poorly known;
    \item the coefficient $\alpha \in \mathbb{R}$ is constant such that the three terms in Eq. \eqref{ul} are of the same magnitude.
\end{itemize}
Despite many successful concrete case studies (see numerous references in \cite{mfc1,mfc2}, and, \textit{e.g.}, \cite{abib,ait,cheng,compar,coskun,ornl,carva,hamon,li,alpha,michel,sancak,forestry,alphavar,he,anae,heg,kenas,yang,wei,pwm,chine} for some recent publications), there are some difficulties:
\begin{itemize}
    \item the determination the order $\nu$ of derivation (see \cite{mfc2});
    \item the determination of the coefficient $\alpha$ and the necessity sometimes to allow time variations (see, \textit{e.g.}, \cite{alphavar,alpha});
    \item the passage to the multivariable case, \textit{i.e.}, to multi-inputs and multi-outputs, where $\alpha$ becomes a matrix (see \cite{serre} for a solution in a concrete case-study).
\end{itemize}

The \emph{HEOL}\footnote{The Breton word \emph{heol} means \emph{sun}.} setting improves and streamlines this approach by taking advantage of the \emph{tangent system}, or \emph{variational system}, associated with the simplified flat system, \textit{i.e.}, the linearized system around a reference trajectory of the simplified flat system.  
When $y$ is a flat output, the tangent linear equation reads:
\begin{equation}\label{lintan}
  \sum_{\rm finite} a_\iota \frac{d^\iota}{dt^\iota} \Delta y  
   = b\Delta u 
\end{equation}
where $a_\iota$, $b$ are possibly time-varying coefficients. Now
\begin{enumerate}
    \item the order $\nu$ of derivation is equal to the least $\iota$, $\iota \neq 0$, such that $a_\iota \not\equiv 0$;
    \item $\alpha = \frac{b}{a_\nu}$: it yields a variable $\alpha$ if $\frac{b}{a_\nu}$ is not constant;
    \item the multivariable case may be dealt with via a diagonal matrix $\alpha$, which is obvious from a theoretical standpoint.
\end{enumerate}
HEOL not only helps to bypass the traditional difficulties of MFC, but also provides a straightforward way to close the loop in flatness-based control. 

It is important to define intrinsically the tangent linear system. Remember moreover that the very concept of flatness was discovered via a crane example \cite{mfc1}, for which the traditional state-variable description fails to hold \cite{crane93}.\footnote{See \cite{bonnabel} for a synthesis of the flatness-based setting in industrial applications of cranes. See also \cite{crane} for the flatness of other crane models.} The algebraic standpoint advocated via differential algebra\footnote{See, \textit{e.g.}, \cite{kolchin}, and \cite{bass} for a general overview of this fascinating domain.} for nonlinear systems \cite{flmr1} and elementary module theory for linear systems \cite{fliess90} permits to do it in a clear cut manner, and, perhaps, in much more precise and elegant way than other techniques (see, \textit{e.g.}, \cite{flmr3,flmr4} for the differential geometry of infinite prolongations). A \emph{tangent linear system} is precisely defined via 
K\"{a}hler differentials \cite{eisenbud,johnson}. It yields the {\em homeostat},\footnote{This word is borrowed from Ashby's remarkable device \cite{homeo}, which is of course related to homeostasis.} which replaces the now classic {\em ultra-local model} \cite{mfc1,mfc2} in model-free control.\footnote{The terminology {\em ultra-local model}, which is now quite popular, becomes irrelevant in this new context, where at least a simplified modeling is available.} In this approach
\begin{itemize}
 \item the feedback design leads to {\em intelligent controllers} which are similar to those in \cite{mfc1,mfc2};
    \item the data-driven estimation techniques, which mimic \cite{algidentif1,algidentif2,sira2}, are based on algebraic manipulations stemming from {\em operational calculus} \cite{yosida}, or Laplace transform;
    \item the whole implementation becomes therefore rather easy. 
\end{itemize}

Our paper is organized as follows. Sect.~\ref{prelim} provides the necessary material for our algebraic viewpoint on linear and nonlinear systems.\footnote{This algebraic point of view is not as popular as it was in Kalman's day~\cite{kalman}. To get a taste of algebra, read Shafarevich's excellent introduction~\cite{shaf}.} Sect.~\ref{homeostasis} is devoted to the homeostat and the associated intelligent controllers. A simple computer experiment is presented and discussed in Sect.~\ref{exper}. See Sect. \ref{conclu} for concluding remarks and hints for future investigations.


\section{Algebraic Preliminaries}
\label{prelim}

\subsection{Linear systems}
\label{linear}

An (\emph{ordinary}) \emph{differential ring} $R$ is a commutative ring equipped with a single \emph{derivation} $\frac{d}{dt} = \dot{} $ such that, $\forall a \in R$, $\frac{da}{dt} = \dot{a} \in R$, and, $\forall a, b \in R$, $\frac{d}{dt} (a + b) = \dot{a} + \dot{b}$, $\frac{d}{dt} (ab) = \dot{a} b + a \dot{b}$. A \emph{constant} is an element $c \in R$ such that $\dot{c} = 0$. An (\emph{ordinary}) \emph{differential field} is an (ordinary) differential ring which is a field.

Let $k$ be a differential field. Write $k[\frac{d}{dt}]$ the ring of linear differential operators $\sum_{\rm finite} a_\ell \frac{d^\ell}{dt^\ell}$, $a_\ell \in k$. This ring is obviously commutative if, and only if, $k$ is a field of constants. In the general noncommutative case, any finitely generated left $k[\frac{d}{dt}]$-module $\mathcal{M}$ satisfies the following property (see, \textit{e.g.}, \cite{cohn}), which is classic in the commutative case, 
\begin{equation}
\label{decomp}
\mathcal{M} = \mathcal{F} \bigoplus \mathcal{T}
\end{equation}
where $\mathcal{F}$ (resp. $\mathcal{T}$) is a \emph{free} (resp. \emph{torsion}) finitely generated left $k[\frac{d}{dt}]$-module. Note that for a finitely generated module $\mathcal{M}$ the following two properties are equivalent
\begin{itemize}
    \item $\mathcal{M}$ is torsion,
    \item $\mathcal{M}$ is finite-dimensional as a $k$-vector space.
\end{itemize}
{\bf Notation.} Write ${\rm span}_{k[\frac{d}{dt}]} (S)$ the submodule of $\mathcal{M}$ spanned by $S \subset \mathcal{M}$.
 A {\em linear system} over the differential ground field $k$ is a finitely generated left $k[\frac{d}{dt}]$-module~$\Lambda$.
A \emph{linear control system} over the ground field $k$ is a finitely generated left $k[\frac{d}{dt}]$-module~$\Lambda$ where
there is a finite set $U = \{u_1, \ldots , u_m\} \subset \Lambda$ of {\em control variables} such that the quotient module $\Lambda/{\rm span}_{k[\frac{d}{dt}]}(U)$ is torsion.
The control variables are said to be {\em independent} if, and only if, ${\rm span}_{k[\frac{d}{dt}]}(U)$ is free of rank~$m$.
System $\Lambda$ is said to be {\em controllable}~\cite{fliess90} if, and only if, $\Lambda$ is a free module. Contrary
to the usual approaches, this definition does not depend on any distinction between system
variables and any state space description. It has been proved nevertheless~\cite{fliess90} that for a standard state-variable representation this module-theoretic definition is equivalent to the classic Kalman's approach.
The set of \emph{output variables} 
$Y = \{y_1, . . . , y_p \} \subset \Lambda$ is a finite subset of the system $\Lambda$. The {\em input-output system} $\Lambda$, with
input $U$ and output $Y$ is said to be {\em observable} \cite{fliess90} if, and only if, $\Lambda = {\rm span}_{k[\frac{d}{dt}]}(U, Y)$, {\it i.e.}, any system
variable is a $k$-linear combination of the control and output variables and their derivatives up
to some finite order. It has been proved \cite{fliess90} that this definition is equivalent to the Kalman definition with a standard state-variable description.

\subsection{Nonlinear systems}

\subsubsection{Differential field extension}

Differential fields are assumed to be of characteristic~$0$. A \emph{differential field extension} $L/K$ is given by two differential fields $K$ and $L$, such that $K \subset L$, and the restriction to~$K$ of the derivation of~$L$ coincides with the derivation of~$K$. For simplicity's sake, $L/K$ is assumed to be finitely generated. Write $K \langle S \rangle$ the differential subfield of $L$ generated over $K$ by $S \subset L$. An element $\xi \in L$ is said to be \emph{differentially algebraic over $K$}, or \emph{differentially $K$-algebraic}, if, and only if, it satisfies an algebraic differential equation over $K$, \textit{i.e.}, there exists a polynomial $\pi [x_0, x_1, \dots, x_\nu]$, $\pi \neq 0$, such that $\pi [\xi, \dot{\xi}, \dots, \xi^{(\nu)} ] = 0$; $\xi$ is said to be \emph{differentially transcendental over $K$}, or \emph{differentially $K$-transcendental}, if, and only if, it is not differentially $K$-algebraic. The extension $L/K$ is said to be \emph{differentially algebraic} (resp. \emph{differentially transcendental}) if, and only if, any (resp. at least one) element in $L$ is (resp. is not) differentially $K$-algebraic. A 
set $\left\{\,\xi_i\, \vert\, i \in I\,\right\}$ is said to be \emph{differentially $K$-algebraically independent}, if, and only if, the set $\left\{x_i^{(\nu)} \vert i \in I, \nu = 0, 1, \dots \right\}$ is algebraically $K$-independent. Such an independent set, which is maximal with respect to inclusion, is called a \emph{differential transcendence basis} of the $L/K$. Two such bases have the same cardinality, \textit{i.e.}, the same number of elements, which is called the \emph{differential transcendence degree} of $L/K$ and is denoted diff tr d° $L/K$. The differential field extension $L/K$ is said to be \emph{purely differentially transcendental} if, and only if, it is generated by a differential transcendence basis. The two following properties are equivalent
\begin{itemize}
    \item diff tr d° $L/K$ $=0$; 
    \item the familiar, \textit{i.e.}, non-differential, transcendence degree of $L/K$ is finite.
\end{itemize} 

\subsubsection{Nonlinear systems and differential flatness}

A \emph{system} is a finitely generated differential field extension $\mathcal{D}/k$. In a \emph{control system}
\begin{itemize}
    \item there is a finite set $U = \{u_1, \dots, u_m\}$ of \emph{control variables},
    \item the extension $\mathcal{D}/k\langle U \rangle$ is differentially algebraic.
\end{itemize}
The control variables are said to be \emph{independent} if, and only if, $k\langle U \rangle / k$ is a purely differential transcendental extension where $U$ is a differential transcendence basis. Introduce a finite set of \emph{output variables} $Y = \{y_1, \dots, y_p \} \subset \mathcal{D}$. This input-output system is said to be \emph{observable} if, and only if, the field extension $\mathcal{D}/k\langle U, Y \rangle$ is algebraic.

The system $\mathcal{D}/k$ is said to be \emph{(differentially)} \emph{flat} if, and only if, the algebraic closure of $\mathcal{D}$ is $k$-isomorphic to the algebraic closure of a purely differentially transcendental extension $k\langle Y \rangle/k$. The components of $Y$ are called \emph{flat}, or \emph{linearizing}, outputs.
Note that any flat system is obviously observable with respect to any flat outputs.

\begin{rema}
\label{non-lin}
Consider the linear system $\Lambda$ of Sec.~\ref{linear} as a $k$-vector space. Let ${\rm Sym}_k \Lambda$ be the \emph{symmetric algebra} (see, \textit{e.g.}, \cite{bourbaki1,shaf}) generated by this $k$-vector space. This integral ring may be endowed with the structure of differential ring. Its field of fractions   
${\rm FracSym}_k \Lambda$ define the differential field extension ${\rm FracSym}_k \Lambda /k$. The correspondence between a basis of the free module~$\Lambda$ and the flat outputs of the flat system ${\rm FracSym}_k \Lambda /k$ demonstrates that a linear system is flat if, and only if, it is controllable. Differential flatness may be viewed as an extension of the familiar Kalman controllability.
\end{rema}

\subsubsection{K\"{a}hler differentials}

\emph{K\"{a}hler differentials} were introduced in commutative algebra and algebraic geometry to mimic some features of differential calculus (see, \textit{e.g.}, \cite{eisenbud,shaf}). They have been extended to differential algebra \cite{johnson}. Consider again a finitely generated differential field extension $L/K$, where $K$ and $L$ are of characteristic $0$. Introduce the (K\"ahler) differential ${\rm d}_{L/K}: L \rightarrow \Omega_{L/K}$ where $\Omega_{L/K}$ is a finitely generated left $L[\frac{d}{dt}]$-module, such that\footnote{$\forall a \in L$,  ${\rm d}_{L/K} a \in \Omega_{L/K}$ should intuitively be viewed as a ``small'' variation of $a$.}
\begin{itemize}
\item $\forall a \in L$, ${\rm d}_{L/K} \dot{a} = \frac{d}{dt}{\rm d}_{L/K} a$; 

\item  $\forall a, b \in L$, ${\rm d}_{L/K} (a + b) = {\rm d}_{L/K} a + {\rm d}_{L/K} b$ and ${\rm d}_{L/K} (ab) = a {\rm d}_{L/K} b + b {\rm d}_{L/K} a$;

\item $\forall c \in K$, ${\rm d}_{L/K} c = 0$.

\end{itemize}
The following properties justify the introduction of K\"{a}hler differentials.
\begin{itemize}
    \item A set $\{\eta_1, \dots, \eta_m\}$ is a differential transcendence basis of $L/K$ if, and only if, $\{\rm d_{L/K} \eta_1, \dots, {\rm d}_{L/K} \eta_m\}$ is a maximal set of $L[\frac{d}{dt}]$-linearly independent elements in $\Omega_{L/K}$. Thus the differential transcendence degree of $L/K$ is equal to the rank of the module $\Omega_{L/K}$.

    \item $L/K$ is differentially algebraic if, and only if, $\Omega_{L/K}$ is torsion. A set $\{\chi_1, \dots, \chi_\nu\}$ is a transcendence basis of $L/K$ if, and only if, $\{{\rm d}_{L/K} \chi_1, \dots, {\rm d}_{L/K} \chi_\nu\}$ is a basis of the $L$-vector space $\Omega_{L/K}$.

    \item $L/K$ is an algebraic extension if, and only if, $\Omega_{L/K} = \{0\}$.
\end{itemize}

The \emph{tangent linear system}, or \emph{variational linear system} associated to the system $\mathcal{D}/k$ is the left $\mathcal{D} [\frac{d}{dt}]$-module $\Omega_{\mathcal{D}/k}$ of K\"{a}hler differentials. If $\mathcal{D}/k$ is flat, $\Omega_{\mathcal{D}/k}$ is obviously free: the tangent system is controllable.
\begin{rema}
Possible singularities of flat systems have been investigated \cite{ollivier1,ollivier2}. The prominent r\^{o}le played by tangent linear systems suggests another way to look which is closer to classic algebraic geometry. Consider, for instance, $\dot{y} = uy$. It is flat and $y$ is a flat output. The tangent linear system reads 
\begin{equation}
\label{sing}
\frac{d}{dt} ({\rm d}_{L/K} y) = u {\rm d}_{L/K} y + y{\rm d}_{L/K} u
\end{equation}
It is degenerated at $y = 0$: the control variable ${\rm d}_{L/K} u$ disappears in Eq. \eqref{sing}. Thus $y = 0$ should be called a singularity. 
\end{rema}

\section{Homeostat}
\label{homeostasis}

\subsection{The monovariable case}

\subsubsection{Preliminary calculations}
Consider a control system $\Sigma$ with a single input (resp. output) variable $u$  (resp. $y$). Assume that it is flat with flat output $y$. This is equivalent to saying that $u$ is algebraic over $k\langle y\rangle$ but not over~$k(y)$. It yields the differential equation
\begin{equation}
\label{E}
E(y, \dot{y}, \dots, y^{(n)}, u) = 0
\end{equation}
where $E$ is a polynomial with coefficients in $k$, where at least one derivative of $y$ appears.
Differentiate~Eq. \eqref{E}: \begin{equation*}
\label{monov}
    \sum_{0 \leqslant \iota \leqslant n} \dfrac{\partial E}{\partial y^{(\iota)}} {\rm d}_{k\langle u, y\rangle/k} y^{(\iota)} + \frac{\partial E}{\partial u} {\rm d}_{k\langle u, y\rangle/k} u = 0
\end{equation*}
Let $\nu$, $0 < \nu \leqslant n$, be the smallest integer such that 
$\dfrac{\partial E}{\partial y^{(\nu)}} \neq 0$. Then
\begin{equation}
\label{interm}
\frac{d^\nu}{dt^\nu} \left({\rm d}_{k\langle u, y\rangle/k} y\right) =  \mathfrak{F} + \mathfrak{a}\,{\rm d}_{k\langle u, y\rangle/k}u   
\end{equation}
where
\begin{eqnarray*}
    \mathfrak{F} &=& - \sum_{\iota \neq \nu} \frac{\frac{\partial E}{\partial y^{(\iota)}}}{\frac{\partial E}{\partial y^{(\nu)}}}{\rm d}_{k\langle u, y\rangle/k} y^{(\iota)} \\
    \mathfrak{a} &=& - \frac{\frac{\partial E}{\partial u}}{\frac{\partial E}{\partial y^{(\nu)}}}
\end{eqnarray*}
The \emph{homeostat}, which is replacing the ultra-local model \cite{mfc1,mfc2}, is deduced from Eq.~\eqref{interm}:
\begin{equation}
\label{driver}
\frac{d^\nu}{dt^\nu} \Delta y  =  F + \alpha \Delta u   
\end{equation}
There
\begin{itemize}
    \item $\Delta y = y - y^\star$, $\Delta u = u - u^\star$, where $y^\star$ is a reference trajectory for the flat system $\Sigma$ and $u^\star$ the \emph{(corresponding) nominal control};
    \item $F = \mathfrak{F} + G$, where $G$ stands for all the mismatches and disturbances.
    \item $\alpha = \mathfrak{a}$ is evaluated on~$y^\star$, and may be time-varying.
\end{itemize}

\subsubsection{Some data-driven calculations}

To encompass systems with time-varying coefficients, let the ground field $k$ be, for instance, a field of meromorphic functions of the variable $t$, such that, $\forall t \in \mathbb{R}$, the coefficients of their Laurent expansions are real. The coefficients of $E$ in Eq. \eqref{E} are, therefore, real-valued functions when they are defined.

In order to estimate $F$ in Eq. \eqref{driver}, we will use, like \cite{mfc1}, classic operational calculus (see, \textit{e.g.}, \cite{yosida}), and the well-known fact that any integrable real-valued function may be approximated by a step function, \textit{i.e.}, a piecewise constant function. Replace Eq. \eqref{driver} by
\begin{equation*}
\label{laplace}
s^\nu \mathcal{Y} - I = \frac{\Phi}{s} + \mathcal{V}
\end{equation*}
where 
\begin{itemize}
    \item $\Phi \in \mathbb{R}$ is a constant to be determined;
    \item $\mathcal{Y}$ (resp. $\mathcal{V}$) is the operational analogue, often called Laplace transform, of $\Delta y$ (resp. $\alpha \Delta u$);
    \item $I \in \mathbb{R}[s]$, is a polynomial of degree less or equal to $\nu - 1$, corresponds to the initial conditions of $\Delta y,\Delta \dot y,\ldots,\Delta y^{(\nu-1)}$.
\end{itemize}
To get rid of $I$, \textit{i.e.}, of the poorly known initial conditions, derive both sides $\nu$ times with respect to $s$, i.e. apply the operator $\frac{d^\nu}{ds^\nu}$. Remember \cite{yosida} that $\frac{d^\nu}{ds^\nu}$ corresponds in the time domain to the multiplication by $(- t)^\nu$. Positive powers of $s$ correspond to time-derivatives. Multiply therefore both sides by $s^{- \mu}$, where $\mu > 0$ is large enough. It yields $\Phi$ as a $\mathbb{R}[\frac{1}{s}]$-linear combination of $\frac{d^\nu \mathcal{V}}{ds^\nu}$, and $\frac{d^\iota \mathcal{Y}}{ds^\iota}$, $0 \leqslant \iota \leqslant \nu$. 

For $\nu = 1$, the operational analogue of Eq.~\eqref{driver} reads
\begin{equation*}
    s\mathcal{Y} - {\Delta}y(0) = \frac{\Phi}s + \mathcal{V}
\end{equation*}
Derive both sides w.r.t. $s$:
\begin{equation*}
    \mathcal{Y} + s\frac{d\mathcal{Y}}{ds} = -\frac{\Phi}{s^2} + \frac{d\mathcal{V}}{ds}
\end{equation*}
Multiply both sides by $s^{-2}$:
\begin{equation*}
-\frac{\Phi}{s^4} = \frac {1}{s^2} \mathcal{Y} + \frac 1 s \frac{d\mathcal{Y}}{ds} - \frac 1{s^2} \frac{d\mathcal{V}}{ds}
\end{equation*}
It yields in the time domain a data-driven real-time estimator $F_{\rm{est}}$:
\begin{equation*}
\label{estim1} F_{\rm est} = - \frac{6}{T^3} \int_{0}^{T} \left( (T - 2 \sigma)\Delta \tilde y(\sigma) {+} \sigma (T - \sigma)\tilde\alpha(\sigma)\Delta \tilde u(\sigma)\right)d\sigma
\end{equation*}
where 
\begin{itemize}
\item the time lapse $T > 0$  is ``small.''
\item $\Delta\tilde y(\sigma)=\Delta y(\sigma+t-T)$, $\tilde\alpha(\sigma)\Delta \tilde u(\sigma)= \alpha(\sigma+t-T)\Delta u(\sigma+t-T)$.
\end{itemize}

For $\nu=2$, analogous calculations give \cite{mfc2}:
\begin{eqnarray*}
    F_{\rm est} &=& \dfrac{60}{T^5}\left[
            \int_0^T \left(\big(T-\sigma\big)^2 -4\big(T-\sigma\big)\sigma + \sigma^2\right)\Delta \tilde y(\sigma)d\sigma
            \right.\nonumber \\
    &&      \mbox{} -
       \left.\dfrac{1}2 \int_0^T (T-\sigma)^2\sigma^2\tilde\alpha(\sigma)\Delta \tilde u(\sigma)d\sigma\right]
  \end{eqnarray*}

\subsubsection{Intelligent controllers}\label{intelcontr}
Introduce \cite{mfc1}, when $\nu = 1$, the \emph{intelligent proportional} controller, or \emph{iP},
\begin{equation}
\label{ip}
\Delta u = - \frac{F_{\text{est}} + K_P \Delta y }{\alpha}    
\end{equation}
where $K_P \in \mathbb{R}$ is the \emph{gain}. Combine Eqs. \eqref{driver} and \eqref{ip}:
$$
\frac{d}{dt} (\Delta y) + K_P \Delta y = F - F_{\text{est}}
$$
If 
\begin{itemize}
    \item the estimate of $F$ is ``good'', \textit{i.e.}, $F - F_{\text{est}} \approx 0$;
    \item $K_P >0$, 
\end{itemize}
then $\displaystyle\lim_{t \to +\infty} \Delta y \approx 0$. This local stability result is easily extended \cite{mfc1,mfc2} to the case $\nu = 2$ via the \emph{intelligent proportional-derivative} controller, or \emph{iPD},
\begin{equation}
\label{ipd}
\Delta u = - \frac{F_{\text{est}} + K_P \Delta y + K_D \frac{d}{dt} (\Delta y)}{\alpha}
\end{equation}
where the gains $K_P, K_D \in \mathbb{R}$ are chosen such that the the roots of $s^2 + K_D s + K_P$ have strictly negative real parts.

\begin{rema}
The extension of \emph{Riachy's trick} \cite{mfc2} to Eq. \eqref{ipd}, which is straightforward, permits to avoid the calculation of the derivative $\frac{d}{dt} (\Delta y)$.
\end{rema}

\subsection{The multivariable case}


Let $\mathcal{D}/k$ a flat multivariable system $\mathcal{D}$ with $m$ independent control variables $U = \{u_1, \dots, u_m\}$ and a flat output $Y = \{y_1, \dots, y_p \}$. Then 
\begin{itemize}
    \item diff tr d° $\mathcal{D}/k = m$, since diff tr d° $\mathcal{D}/k\langle U\rangle = 0$,
    \item diff tr d° $\mathcal{D}/k =$ diff tr d° $k\langle Y\rangle/k = p$.
\end{itemize}
Thus $p=m$, the number of flat output variables is equal to the number of independent control variables.

Every component of $U$ is algebraic over $k\langle Y\rangle$  but not over $k(Y)$. Therefore there exists differential equations of the form:

Assume that the components of $U$ are algebraic over $k\langle Y\rangle$  but not over $k(Y)$.\footnote{If not it leads to algebraic equations and not to differential ones.} Therefore there exists differential equations of the form:
\begin{equation*}
    E_j(Y,\dot Y,\dots,Y^{(\nu_j)}, u_j ) =0,\qquad j=1,\ldots,m
\end{equation*}
where $E_j$ is a polynomial with coefficients in~$k$ in which at least a derivative of one component of~$y$ appears. Taking now the K\"ahler differential of the $E_j$'s and, up to a renumbering of the components of the flat output, one obtains:
\begin{equation*}
    \forall j=1,\ldots,m,\qquad \sum_{1\leqslant l\leqslant m}\quad\sum_{0 \leqslant \iota_l \leqslant \nu_j} \frac{\partial E_j}{\partial y_l^{(\iota_l)}} {\rm d}_{k\langle U, Y\rangle/k} y_l^{(\iota_l)} + \frac{\partial E_j}{\partial u_j} {\rm d}_{k\langle U, Y\rangle/k} u_j = 0
\end{equation*}
Let $\mu_j$, $0 < \mu_j \leqslant \nu_j$, be the smallest integer such that 
$\dfrac{\partial E_j}{\partial y_j^{(\mu_j)}} \neq 0$. Then
\begin{equation}
\label{interm:multi}
\frac{d^{\mu_j}}{dt^{\mu_j}} \left({\rm d}_{k\langle U, Y\rangle/k} y_j\right) =  \mathfrak{F}_j + \mathfrak{a}_j\, {\rm d}_{k\langle U, Y\rangle/k}u_j   
\end{equation}
where
\begin{eqnarray*}
\mathfrak{F}_j &=& - \displaystyle\sum_{1\leqslant l\leqslant m}\ \sum_{\iota_l \neq \mu_j} \frac{\frac{\partial E_j}{\partial y_l^{(\iota_l)}}}{\frac{\partial E_j}{\partial y_l^{(\mu_j)}}}\\
\mathfrak{a}_j &=& - \frac{\frac{\partial E_j}{\partial u_j}}{\frac{\partial E}{\partial y_j^{({\mu_j})}}}
\end{eqnarray*}
The \emph{homeostat}, which is deduced from Eq.~\eqref{interm:multi}, reads
\begin{equation*}
\label{driver1}
{\frac{d^{\mu_j}}{dt^{\mu_j}}}\Delta y_j  =  F_j + \alpha_j \Delta u_j,\qquad j=1,\ldots,m
\end{equation*}
with
\begin{itemize}
    \item $\Delta y_j = y_j - y^\star_j$, $\Delta u_j = u_j - u^\star_j$, where $Y^\star = \{y_1^\star,\ldots,y_m^\star\}$ is a reference trajectory for the flat system $\mathcal{D}$ and $U^\star = \{u_1^\star,\ldots,u_m^\star\}$ the \emph{(corresponding) nominal control};
    \item $F_j = \mathfrak{F}_j + G_j$, where $G_j$ stand for all the mismatches and disturbances.
\end{itemize}
The extension of Section~\ref{intelcontr} to the multivariable case is straightforward. 

\section{A computer experiment}
\label{exper}

Consider the flat system
$$
\begin{cases}
\dot x_1=x_1+x_1^2u_1\\
\dot x_2=x_3\\
\dot x_3=x_4\\
\dot x_4=-x_4+x_3+x_2+x_1u_1u_2\\
y_1=x_1\\
y_2=x_2
\end{cases}
$$
where $y_1$, $y_2$ are flat outputs. 
The nominal control variables are given by
$$u_1^\star=\frac{\dot y_1^\star-y_1^\star}{{y_1^\star}^2}$$
and
$$u_2^\star=\frac{\dddot y_2^\star+\ddot y_2^\star-\dot y_2^\star-y_2^\star}{y_1^\star u_1^\star}$$
The homeostat becomes
$$
\begin{cases}
\frac{d}{dt}(\Delta y_1) =F_1+{y_1^\star}^2\Delta u_1\\
\frac{d^2}{dt^2}(\Delta y_1)=F_2+\left(\frac{\dot y_1^\star}{y_1^\star}-1\right)\Delta u_2
\end{cases}
$$

\begin{figure*}[!ht]
\centering%
\subfigure[\footnotesize Control 1 (blue --) and nominal control 1 (red - -) ]
{\epsfig{figure=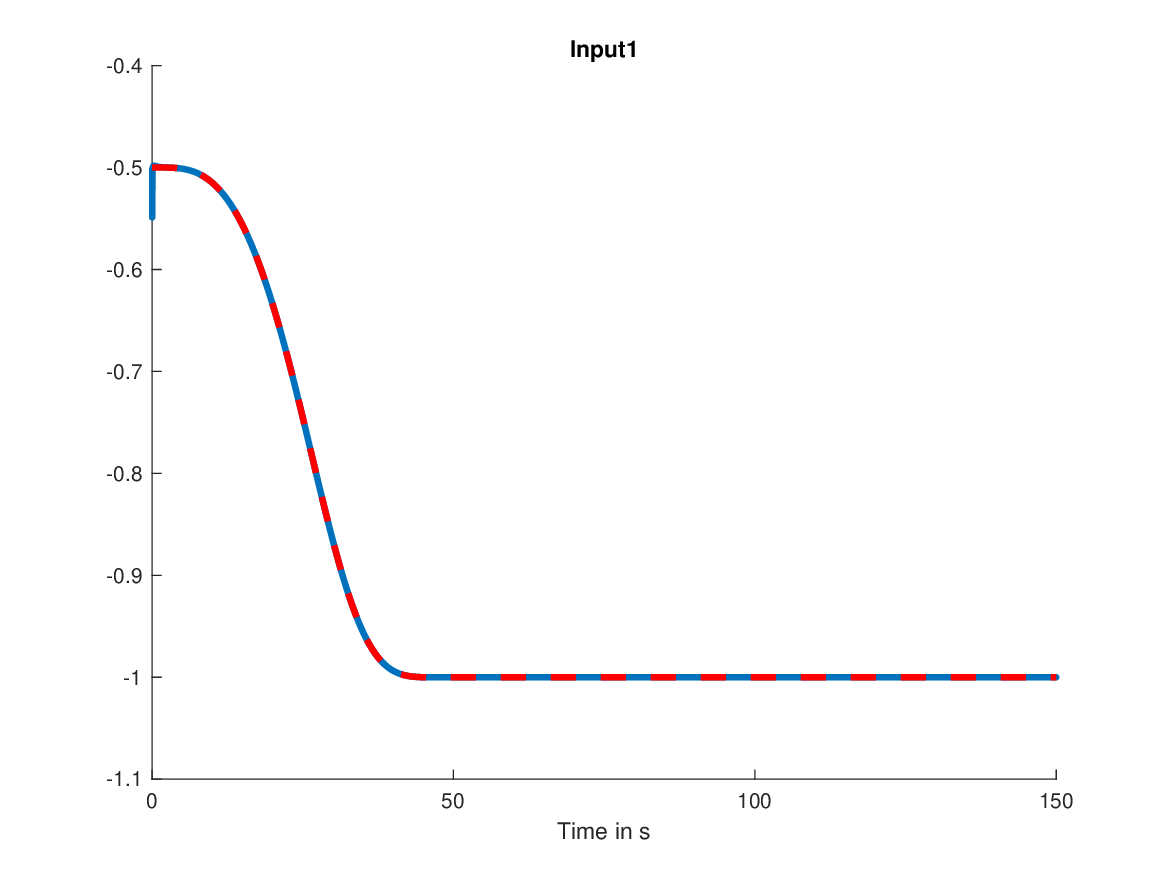,width=0.48\textwidth}}
\subfigure[\footnotesize Control 2  (blue --) and nominal control 2(red - -) ]
{\epsfig{figure=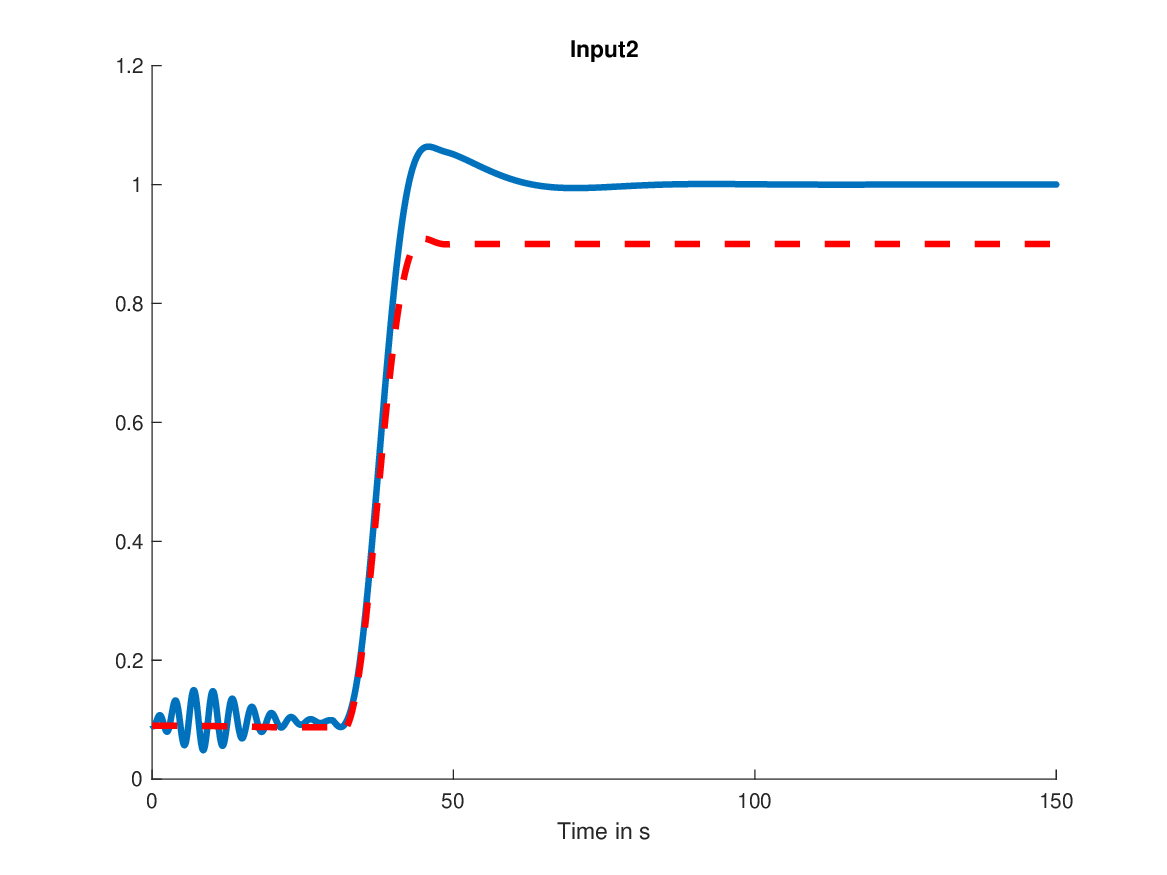,width=0.48\textwidth}}
\\
\subfigure[\footnotesize Output 1 (blue --) and trajectory 1 (red - -)]
{\epsfig{figure=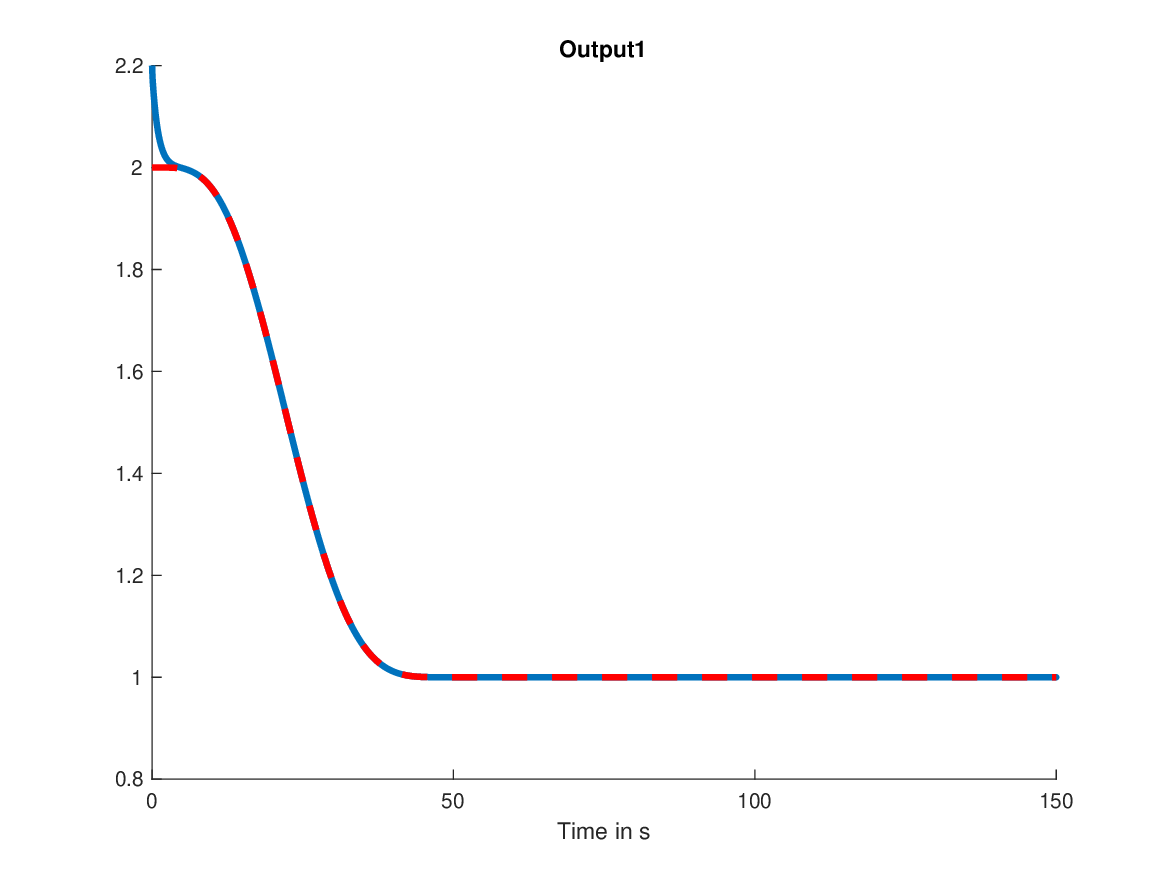,width=0.48\textwidth}}
\subfigure[\footnotesize Output 2 (blue --) and trajectory 2 (red - -)]
{\epsfig{figure=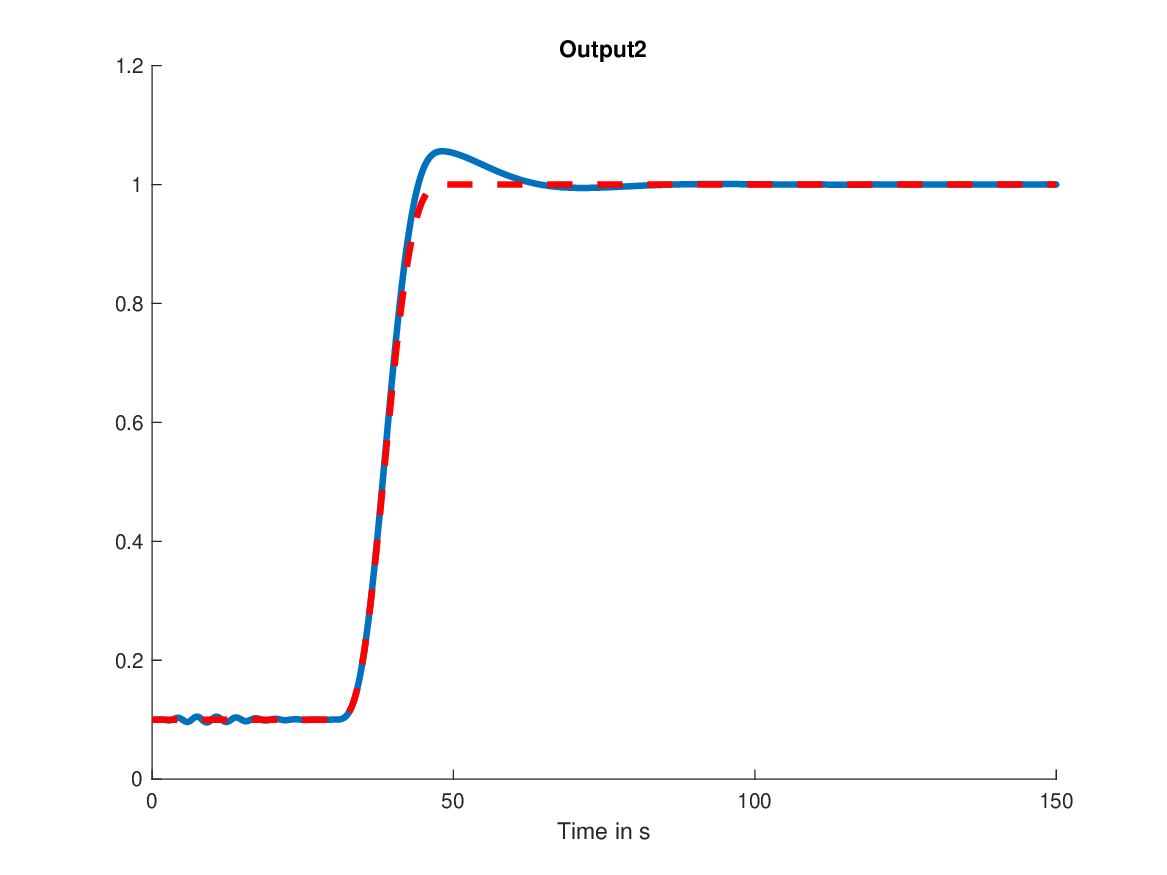,width=0.48\textwidth}}

\caption{Control evaluation}\label{RM2S2}
\end{figure*}

\noindent Close the loop for the first (resp. second) equation such that the root (resp. double root) of the characteristic polynomial is $-1$ (resp. $-0.15$). The simulation duration is $150$\,s. The sampling period is $10$\,ms. The following mismatches are introduced to show the robustness of our control strategy:
\begin{itemize}
\item 
$y_1(0)= 1.1 y_1^\star(0)$, $y_2(0)= y_2^\star(0)$;
\item  $u_2^\star=\frac{\dddot y_2^\star+\ddot y_2^\star-1.1\dot y_2^\star-0.9y_2^\star}{y_1^\star u_1^\star}$
\end{itemize}
The results displayed in Figure \ref{RM2S2} are quite satisfactory. 




\section{Conclusion}
\label{conclu}

The HEOL setting, \textit{i.e.}, the introduction of homeostats deduced from the tangent linear system, suggests elementary solutions to questions, which are crucial from a practical viewpoint and have been around for many years, like  feedback control of flat systems and implementation issues in model-free control. There are of course other issues which might benefit from our approach. It has been observed (see, \textit{e.g.}, \cite{optim} for a recent contribution) that flatness-based control greatly simplifies optimal control. The HEOL combination should bring further improvements (see, \textit{e.g.}, \cite{euler}). Convincing concrete illustrations should be available soon.

Let's conclude with some more general considerations. In applied sciences, too, an appropriate formalism might be of paramount importance in order to trivialize technical investigations that seemed before beyond the reach. This is the aim of the present paper in control engineering.

\bibliographystyle{crplain}


\bibliography{flat-rev-cr-math}\end{document}